\documentclass[12pt,draft]{amsart}
\usepackage[cp1251]{inputenc}
\usepackage[T2A]{fontenc}
\usepackage[english]{babel}
\usepackage{amsmath,amsfonts,amssymb}
\usepackage{geometry}

\newtheorem{proposition}{Proposition}

\theoremstyle{remark}

\theoremstyle{definition}

\begin{document}

\title[Parameter-elliptic problems and interpolation]
{Parameter-elliptic problems and\\interpolation with a function parameter}


\author[A. Anop]{Anna V. Anop}

\address{Institute of Mathematics, National Academy of Sciences of Ukraine,
3 Tereshchenkivs'ka, Kyiv, 01601, Ukraine}

\email{ahlv@ukr.net}


\author[A. Murach]{Aleksandr A. Murach}

\address{Institute of Mathematics, National Academy of Sciences of Ukraine,
3 Tereshchenkivs'ka, Kyiv, 01601, Ukraine}

\email{murach@imath.kiev.ua}






\begin{abstract}
Parameter-elliptic boundary-value problems are investigated on the extended Sobolev
scale. This scale consists of all Hilbert spaces that are interpolation spaces with
respect to the Hilbert Sobolev scale. The latter are the H\"ormander spaces
$B_{2,k}$ for which the smoothness index $k$ is an arbitrary radial function
RO-varying at~$+\infty$. We prove that the operator corresponding to this problem
sets isomorphisms between appropriate H\"ormander spaces  provided that the absolute
value of the parameter is large enough. For solutions to the problem, we establish
two-sided estimates, in which the constants are independent of the parameter.
\end{abstract}

\maketitle

\section{Introduction}\label{sec1}

Interpolation between normed spaces is an efficient and convenient method in the
theory of operators. This is due to the fact that boundedness of linear operators,
their isomorphism and Fredholm properties remain valid under the interpolation
between the spaces on which these operators act. In the theory of differential
operators, various methods of the interpolation are used to prove theorems about
collections of isomorphisms or Fredholm mappings realized by differential operators
on spaces that are intermediate for a given couple of Sobolev spaces. This approach
is taken in the monographs by Yu.~M.~Berezansky \cite[Chap.~III,
Sec.~6]{Berezansky68}, J.-L.~Lions and E.~Magenes~\cite{LionsMagenes72i,
LionsMagenes72ii}, H.~Triebel \cite{Triebel95}, V.~A.~Mikhailets and
A.~A.~Murach~\cite{MikhailetsMurach10, MikhailetsMurach14}.

In the present paper, we give an application of interpolation with a function
parameter between Hilbert spaces to an important class of parameter-elliptic
boundary-value problems. This class was selected and investigated on the Sobolev
scale by S.~Agmon, L.~Nirenberg \cite{Agmon62, AgmonNirenberg63} and
M.~S.~Agranovich, M.~I.~Vishik \cite{AgranovichVishik64}. They proved that the class
has the following fundamental property. If the complex parameter is large enough in
absolute value, then the operator corresponding to the parameter-elliptic problem
sets an isomorphism between appropriate Sobolev spaces. Moreover, the solutions to
this problem admit a two-sided estimate with constants that do not depend on the
parameter. These properties have important applications to the spectral theory of
elliptic operators and parabolic differential equations. Various and more general
classes of parameter-elliptic problems are investigated by M.~S.~Agranovich
\cite{Agranovich90, Agranovich92}, G.~Grubb \cite[Ch. 2]{Grubb96}, A.~N.~Kozhevnikov
\cite{Kozhevnikov96, Kozhevnikov97, Kozhevnikov01}, R.~Denk, R.~Mennicken, and
L.~R.~Volevich \cite{DenkMennickenVolevich98, DenkMennickenVolevich01}, R.~Denk and
M.~Fairman \cite{DenkFairman10} and others (see also survey \cite{Agranovich97}).

The purpose of this paper is to prove the above-mentioned fundamental property for
the class of all Hilbert spaces that are interpolation spaces with respect to the
Hilbert Sobolev scale. It follows from V.~I.~Ovchinnikov's theorem \cite[Sec.
11.4]{Ovchinnikov84} that these spaces are obtained by interpolation with a function
parameter between inner product Sobolev spaces. The class of these interpolation
spaces is constructively described in \cite[Sec.~2.4.2]{MikhailetsMurach10,
MikhailetsMurach14}. It consists of the inner product H\"ormander spaces $B_{2,k}$
\cite[Sec.~2.2]{Hermander63}  for which the smoothness index $k$ is an arbitrary
radial function RO-varying at $+\infty$. This class is naturally to call the
extended Sobolev scale (by means of interpolation spaces).

Note that various classes of elliptic differential operators given on manifolds
without boundary are investigated on this scale in \cite{09UMJ3, 13MFAT1, 12UMJ11}
and \cite[Sec.~2.4.3]{MikhailetsMurach10, MikhailetsMurach14}.

\section{Statement of the problem}\label{sec2}

Let $\Omega$ be a bounded domain in $\mathbb{R}^{n}$, with $n\geq2$. Suppose that
its boundary $\Gamma:=\partial\Omega$ is an infinitely smooth closed manifold of
dimension $n-1$. Let $\nu(x)$ denote the unit vector of the inner normal to
$\partial\Omega$ at a point $x\in\Gamma$. As usual,
$\overline{\Omega}:=\Omega\cup\Gamma$.

Let $q\geq1$ and $m_{1},\ldots,m_{q}\geq0$ be arbitrarily chosen integers. We
consider the boun\-dary-value problem
\begin{equation}\label{4f1}
A(\lambda)\,u=f\quad\mbox{in}\quad\Omega,\quad\quad
B_{j}(\lambda)\,u=g_{j}\quad\mbox{on}\quad\Gamma,\quad j=1,\ldots,q,
\end{equation}
that depends on the parameter $\lambda\in\mathbb{C}$ as follows:
\begin{equation}\label{4f2}
A(\lambda):=\sum_{r=0}^{2q}\,\lambda^{2q-r}A_{r},\quad\quad
B_{j}(\lambda):=\sum_{r=0}^{m_{j}}\,\lambda^{m_{j}-r}B_{j,r}.
\end{equation}
Here, all $A_{r}$ and $B_{j,r}$ are linear partial differential expressions whose
orders do not exceed~$r$. We write these expressions in the form
\begin{gather}\label{4f3}
A_{r}:=A_{r}(x,D):=\sum_{|\mu|\leq r}a_{r,\mu}(x)D^{\mu},\quad x\in\overline{\Omega},\\
B_{j,r}:=B_{j,r}(x,D):=\sum_{|\mu|\leq r}b_{j,r,\mu}(x)D^{\mu},\quad x\in\Gamma.
\label{4f4}
\end{gather}
Their coefficients are complex-valued infinitely smooth functions; i.e., all
$a_{r,\mu}\in C^{\infty}(\overline{\Omega})$ and $b_{j,r,\mu}\in
C^{\infty}(\Gamma)$. Put $B(\lambda):=(B_{1}(\lambda),\ldots,B_q(\lambda))$.

Note that we use the standard notation in \eqref{4f3}, \eqref{4f4}, and below;
namely, for the multi-index $\mu=(\mu_{1},\ldots,\mu_{n})$ we let
$|\mu|:=\mu_{1}+\ldots+\mu_{n}$ and $D^{\mu}:=D_{1}^{\mu_{1}}\ldots
D_{n}^{\mu_{n}}$, with $D_{k}:=i\,\partial/\partial x_{k}$ for $k\in\{1,\ldots,n\}$
and $x=(x_{1},\ldots,x_{n})\in\mathbb{R}^{n}$. Moreover, we put
$\xi^{\mu}:=\xi_{1}^{\mu_{1}}\ldots\xi_{n}^{\mu_{n}}$ for the frequency variables
$\xi=(\xi_{1},\ldots,\xi_{n})\in\mathbb{C}^{n}$, which are dual to the spatial
variables $x=(x_{1},\ldots,x_{n})$ with respect to the Fourier transform
$$
\widehat{w}(\xi):=(\mathcal{F}w)(\xi):=\frac{1}{(2\pi)^{n/2}}
\int\limits_{\mathbb{R}^{n}}e^{i\xi x}\,w(x)\,dx.
$$

Following M.~S.~Agranovich and M.~I.~Vishik \cite[Sec. 4]{AgranovichVishik64} (see
also survey \cite[Sec. 3.1]{Agranovich97}), we recall the notion of
parameter-ellipticity in connection with the boundary-value problem \eqref{4f1}.

Let us associate certain homogeneous polynomials in
$(\xi,\lambda)\in\mathbb{C}^{n+1}$ with partial differential expressions
\eqref{4f2}. Namely, we set
\begin{gather*}
A^{(0)}(x,\xi,\lambda):=
\sum_{r=0}^{2q}\,\lambda^{2q-r}\sum_{|\mu|=r}a_{r,\mu}(x)\,\xi^{\mu},
\quad\mbox{with}\quad x\in\overline{\Omega},\\
B^{(0)}_{j}(x,\xi,\lambda):=
\sum_{r=0}^{m_{j}}\,\lambda^{m_{j}-r}\sum_{|\mu|=r}b_{j,r,\mu}(x)\,\xi^{\mu},
\quad\mbox{with}\quad x\in\Gamma.
\end{gather*}

Let $K$ be a fixed closed angle on the complex plain with vertex at the origin; this
angle may degenerate into a ray.

The boundary-value problem \eqref{4f1} is called parameter-elliptic in the angle $K$
if the following two conditions are satisfied:
\begin{enumerate}
\item[i)] $A^{(0)}(x,\xi,\lambda)\neq0$ for each $x\in\overline{\Omega}$ and all
$\xi\in\mathbb{R}^{n}$ and $\lambda\in K$ with $|\xi|+|\lambda|\neq0$.
\item[ii)] Let $x\in\Gamma$, $\xi\in\mathbb{R}^{n}$, and $\lambda\in K$
be arbitrarily chosen so that $\xi$ is tangent to $\Gamma$ at $x$ and that
$|\xi|+|\lambda|\neq0$. Then the polynomials $B^{(0)}_{j}(x,\xi+\tau\nu(x),\lambda)$
in~$\tau$, $j=1,\ldots,q$, are linearly independent modulo
$\prod_{j=1}^{q}(\tau-\tau^{+}_{j}(x,\xi,\lambda))$. Here,
$\tau^{+}_{1}(x,\xi,\lambda),\ldots,\tau^{+}_{q}(x,\xi,\lambda)$ are all
$\tau$-roots of $A^{(0)}(x,\xi+\tau\nu(x),\lambda)$ with
$\mathrm{Im}\,\tau>\nobreak0$, each root being taken the number of times equal to
its multiplicity.
\end{enumerate}

In connection with Condition ii), it is relevant to mention the following fact
\cite[Proposition 2.2]{AgranovichVishik64}. If $A^{(0)}(x,\xi,\lambda)$ satisfies
Condition i), then the polynomial $A^{(0)}(x,\xi+\tau\nu(x),\lambda)$ in $\tau$ has
$q$ roots with $\mathrm{Im}\,\tau>0$ and $q$ roots with $\mathrm{Im}\,\tau<0$ (if we
take their multiplicity into account). Thus, the above definition of
parameter-ellipticity is reasonably formulated.

Some examples of parameter-elliptic boundary-value problems are given in \cite[Sec.
3.1 b]{Agranovich97}. For instance, the following boundary-value problem
$$
\Delta u+\lambda^{2}u=f\;\;\mbox{in}\;\;\Omega,\quad\quad\frac{\partial
u}{\partial\nu}-\lambda u=g\;\;\mbox{on}\;\;\Gamma
$$
is parameter-elliptic in the angle $K$ provided that $K$ contains neither
$\mathbb{R}_{+}$ nor $\mathbb{R}_{-}$, where
$\mathbb{R}_{\pm}:=\{\lambda\in\mathbb{R}:\lambda\gtrless0\}$. Note also that, if
$A(\lambda)$ satisfies Condition i), then the Dirichlet boundary-value problem for
the partial differential equation $A(\lambda)=f$ is parameter-elliptic in the angle
$K$.

We investigate properties of the parameter-elliptic boundary-value problem
\eqref{4f1} on the extended Sobolev scale.

\section{The extended Sobolev scale}\label{sec3}

This scale \cite{13UMJ3}, \cite[Sec. 2.4.2]{MikhailetsMurach14} consists of the
inner product isotropic H\"ormander spaces $H^{\varphi}$ for which an arbitrary
function parameter $\varphi\in\mathrm{RO}$ serves as a smoothness index. Here, RO is
the class of all Borel measurable functions
$\varphi:[1,\infty)\rightarrow(0,\infty)$ such that
\begin{equation}\label{4f5}
c^{-1}\leq\frac{\varphi(\lambda t)}{\varphi(t)}\leq c\quad\mbox{for arbitrary}\quad
t\geq1\quad\mbox{and}\quad\lambda\in[1,a]
\end{equation}
with some numbers $a=a(\varphi)>1$ and $c=c(\varphi)\geq1$ (certainly, $a$ and $c$
are independent of both $t$ and $\lambda$). Such functions are said to be RO-varying
in the sense of V.~G.~Avakumovi\'c \cite{Avakumovic36} and are sufficiently
investigated \cite{BinghamGoldieTeugels89, Seneta76}.

The class RO admits the simple description
\begin{equation*}
\varphi\in\mathrm{RO}\quad\Leftrightarrow\quad\varphi(t)=
\exp\biggl(\beta(t)+\int_{1}^{t}\frac{\gamma(\tau)}{\tau}\;d\tau\biggr), \;\;t\geq1,
\end{equation*}
where the real-valued functions $\beta$ and $\gamma$ are Borel measurable and
bounded on $[1,\infty)$. Note also that condition \eqref{4f5} is equivalent to the
bilateral inequality
\begin{equation}\label{4f6}
c_{0}\lambda^{s_{0}}\leq\frac{\varphi(\lambda t)}{\varphi(t)}\leq
c_{1}\lambda^{s_{1}}\quad\mbox{for arbitrary}\quad
t\geq1\quad\mbox{and}\quad\lambda\geq1,
\end{equation}
where $c_{0}$ and $c_{1}$ are certain positive numbers, which do not depend on $t$
and $\lambda$. Hence, for every function $\varphi\in\mathrm{RO}$, we can define the
lower and the upper Matuszewska indices \cite{Matuszewska64} by the formulas
\begin{gather*}
\sigma_{0}(\varphi):=\sup\{s_{0}\in\mathbb{R}:\,\mbox{the left-hand inequality in
\eqref{4f6} holds}\},\\
\sigma_{1}(\varphi):=\inf\{s_{1}\in\mathbb{R}:\,\mbox{the right-hand inequality in
\eqref{4f6} holds}\},
\end{gather*}
with $-\infty<\sigma_{0}(\varphi)\leq\sigma_{1}(\varphi)<\infty$ (see \cite[Theorem
2.2.2]{BinghamGoldieTeugels89}).

Now let $\varphi\in\mathrm{RO}$ and introduce the function spaces $H^{\varphi}$ over
$\mathbb{R}^{n}$ and then over $\Omega$ and~$\Gamma$.

By definition, the linear space $H^{\varphi}(\mathbb{R}^{n})$ consists of all
distributions $w\in\mathcal{S}'(\mathbb{R}^{n})$ such that their Fourier transform
$\widehat{w}$ is locally Lebesgue integrable over $\mathbb{R}^{n}$ and satisfies the
condition
$$
\int_{\mathbb{R}^{n}}\varphi^2(\langle\xi\rangle)\,|\widehat{w}(\xi)|^2\,d\xi
<\infty.
$$
Here, as usual, $\mathcal{S}'(\mathbb{R}^{n})$ is the complex linear topological
space of tempered distributions given in $\mathbb{R}^{n}$, and
$\langle\xi\rangle:=(1+|\xi|^{2})^{1/2}$ is the smoothed modulus of
$\xi\in\mathbb{R}^{n}$. The inner product in $H^{\varphi}(\mathbb{R}^{n})$ is
defined as follows
$$
(w_1,w_2)_{H^{\varphi}(\mathbb{R}^{n})}:=
\int_{\mathbb{R}^{n}}\varphi^2(\langle\xi\rangle)\,
\widehat{w_1}(\xi)\,\overline{\widehat{w_2}(\xi)}\,d\xi,
$$
with $w_1,w_2\in H^{\varphi}(\mathbb{R}^{n})$. It induces the norm
$\|w\|_{H^{\varphi}(\mathbb{R}^{n})}:=(w,w)_{H^{\varphi}(\mathbb{R}^{n})}^{1/2}$.

The space $H^{\varphi}(\mathbb{R}^{n})$ is a Hilbert and isotropic case of the
spaces $B_{p,k}$ introduced and systematically investigated by L.~H\"ormander
\cite[Sec. 2.2]{Hermander63} (also see \cite[Sec. 10.1]{Hermander83}). Namely,
$H^{\varphi}(\mathbb{R}^{n})=B_{p,k}$ provided that $p=2$ and
$k(\xi)=\varphi(\langle\xi\rangle)$ for all $\xi\in\mathbb{R}^{n}$. Not that, in the
Hilbert case of $p=2$, the H\"ormander spaces coincide with the spaces introduced
and studied by L.~R.~Volevich and B.~P.~Paneah \cite[Sec.~2]{VolevichPaneah65}.

If $\varphi(t)\equiv t^{s}$, then
$H^{\varphi}(\mathbb{R}^{n})=:H^{(s)}(\mathbb{R}^{n})$ is the inner product Sobolev
space of order $s\in\mathbb{R}$. Generally,
\begin{equation}\label{4f7}
s_{0}<\sigma_{0}(\varphi)\leq\sigma_{1}(\varphi)<s_{1}\;\;\Rightarrow\;\;
H^{(s_1)}(\mathbb{R}^{n})\hookrightarrow H^{\varphi}(\mathbb{R}^{n})\hookrightarrow
H^{(s_0)}(\mathbb{R}^{n}),
\end{equation}
both embeddings being continuous and dense.

Following \cite{13UMJ3}, we call the class of Hilbert function spaces
\begin{equation}\label{4f8}
\bigl\{H^{\varphi}(\mathbb{R}^{n}):\varphi\in\mathrm{RO}\bigr\}
\end{equation}
the extended Sobolev scale over $\mathbb{R}^{n}$.

Its analogs for the Euclidean domain $\Omega$ and for the closed manifold $\Gamma$
are introduced in the standard way on the basis of the class \eqref{4f8} (see
\cite[Sec.~2]{arXiv:1106.2049} and \cite[Sec.~2.4.2]{MikhailetsMurach10,
MikhailetsMurach14}). Let us give the necessary definitions; as above,
$\varphi\in\mathrm{RO}$.

By definition,
\begin{gather}\notag
H^{\varphi}(\Omega):=\bigl\{w\!\upharpoonright\!\Omega:
\,w\in H^{\varphi}(\mathbb{R}^{n})\bigr\},\\
\|u\|_{H^{\varphi}(\Omega)}:=\inf\bigl\{\,\|w\|_{H^{\varphi}(\mathbb{R}^{n})}:\,
w\in H^{\varphi}(\mathbb{R}^{n}),\;w=u\;\,\mbox{in}\;\,\Omega\bigr\}, \label{4f9}
\end{gather}
with $u\in H^{\varphi}(\Omega)$. The linear space $H^{\varphi}(\Omega)$ is Hilbert
and separable with respect to the norm \eqref{4f9} because $H^{\varphi}(\Omega)$ is
the factor space of the separable Hilbert space $H^\varphi(\mathbb{R}^{n})$ by its
subspace
\begin{equation*}
\bigl\{w\in H^{\varphi}(\mathbb{R}^{n}):\,
\mathrm{supp}\,w\subseteq\mathbb{R}^{n}\setminus\Omega\bigr\}.
\end{equation*}
Note that $C^\infty(\overline{\Omega})$ is dense in $H^\varphi(\Omega)$.

By definition, the linear space $H^{\varphi}(\Gamma)$ consists of all distributions
on $\Gamma$ that belong to $H^{\varphi}(\mathbb{R}^{n-1})$ in local coordinates on
$\Gamma$. Namely, arbitrarily chose a finite collection of local charts
$\alpha_{j}:\mathbb{R}^{n-1}\leftrightarrow\Gamma_{j}$, $j=1,\ldots,\varkappa$,
where the open sets $\Gamma_{j}$ form a covering of $\Gamma$. Also choose functions
$\chi_{j}\in C^{\infty}(\Gamma)$, $j=1,\ldots,\varkappa$, that satisfy the condition
$\mathrm{supp}\,\chi_{j}\subset\Gamma_{j}$ and that form a partition of unity on
$\Gamma$. Then
\begin{equation*}
H^{\varphi}(\Gamma):=\bigl\{h\in\mathcal{D}'(\Gamma):\,
(\chi_{j}h)\circ\alpha_{j}\in H^{\varphi}(\mathbb{R}^{n-1})\;\;\mbox{for
every}\;\;j\in\{1,\ldots,\varkappa\}\bigr\}.
\end{equation*}
Here, as usual, $\mathcal{D}'(\Gamma)$ is the complex linear topological space of
all distributions given on~$\Gamma$, and $(\chi_{j}h)\circ\alpha_{j}$ is the
representation of the distribution $\chi_{j}h$ in the local chart $\alpha_{j}$. The
inner product in $H^{\varphi}(\Gamma)$ is defined by the formula
$$
(h_{1},h_{2})_{H^{\varphi}(\Gamma)}:=
\sum_{j=1}^{\varkappa}\,((\chi_{j}h_{1})\circ\alpha_{j},
(\chi_{j}\,h_{2})\circ\alpha_{j})_{H^{\varphi}(\mathbb{R}^{n-1})},
$$
with $h_{1},h_{2}\in H^{\varphi}(\Gamma)$, and induces the norm
$\|h\|_{H^{\varphi}(\Gamma)}:=(h,h)_{H^{\varphi}(\Gamma)}^{1/2}$. The space
$H^\varphi (\Gamma)$ is Hilbert and separable and does not depend (up to equivalence
of norms) on our choice of local charts and partition of unity on~$\Gamma$
\cite[Theorem 2.21]{MikhailetsMurach10, MikhailetsMurach14}. Note also that
$C^{\infty}(\Gamma)$ is dense in $H^{\varphi}(\Gamma)$.

The above-defined function spaces form the extended Sobolev scales
\begin{equation}\label{4f10}
\bigl\{H^{\varphi}(\Omega):\varphi\in\mathrm{RO}\bigr\}\quad\mbox{and}\quad
\bigl\{H^{\varphi}(\Gamma):\varphi\in\mathrm{RO}\bigr\}
\end{equation}
over $\Omega$ and $\Gamma$ respectively. They contain the scales of inner product
Sobolev spaces; namely, if $\varphi(t)\equiv t^s$ for some $s\in\mathbb{R}$, then
$H^\varphi(\Omega)=:H^{(s)}(\Omega)$ and $H^\varphi(\Gamma)=:H^{(s)}(\Gamma)$ are
the Sobolev spaces of order~$s$. Property \eqref{4f7} remains true (with embeddings
being continuous and dense) provided that we replace $\mathbb{R}^{n}$ by $\Omega$ or
$\Gamma$ therein.

\section{The main result}\label{sec4}

Let us state our main result regarding the isomorphism property of the
parameter-elliptic problem \eqref{4f1} considered on the extended Sobolev scale.

To avoid mentioning the argument $t$ in the smoothness indices, we refer to
$\varrho$ as the function $\varrho(t):=t$ of $t\geq1$. Note that, if
$\varphi\in\mathrm{RO}$ and $s\in\mathbb{R}$, then the function $\varrho^s\varphi$
belongs to $\mathrm{RO}$, and its Matuszewska indices satisfy the relation
$\sigma_j(\varrho^s\varphi)=s+\sigma_j(\varphi)$ for each $j\in\{0,1\}$.

The mapping $u\mapsto(A(\lambda)u,B(\lambda)u)$, with $u\in
C^\infty(\overline{\Omega})$, extends uniquely (by continuity) to the bounded linear
operator
\begin{equation}\label{4f11}
(A(\lambda),B(\lambda)):\,H^{\varphi\varrho^{2q}}(\Omega)\rightarrow
H^{\varphi}(\Omega)\oplus\bigoplus_{j=1}^{q}H^{\varphi\varrho^{2q-m_j-1/2}}(\Gamma)
=:\mathcal{H}^\varphi(\Omega,\Gamma)
\end{equation}
for each $\lambda\in\mathbb{C}$ and every function parameter $\varphi\in\mathrm{RO}$
that meets the condition
\begin{equation}\label{4f12}
\sigma_{0}(\varphi)>l:=\max\bigl\{0,\,m_{1}-2q+1/2,\ldots,m_{q}-2q+1/2\bigr\}.
\end{equation}

We need to introduce certain norms, which depend on the parameter
$p:=|\lambda|\geq1$. Let $\alpha\in\mathrm{RO}$, with $\sigma_{0}(\alpha)>0$, and
let $G\in\{\mathbb{R}^{n},\Omega,\Gamma\}$. We define an equivalent Hilbert norm in
the space $H^{\alpha}(G)$ by the formula
\begin{equation}\label{4f13}
\|w\|_{\alpha,p,G}:=\bigl(\,\|w\|_{H^{\alpha}(G)}^{2}+
\alpha^{2}(p)\,\|w\|_{H^{(0)}(G)}^{2}\bigr)^{1/2},
\end{equation}
with $w\in H^{\alpha}(G)$. The equivalence of norms follows from the continuous
embedding $H^{\alpha}(G)\hookrightarrow H^{(0)}(G)$. Note that the space
$H^{(0)}(G)=L_{2}(G)$ consist of all the functions that are square integrable over
$G$.

\medskip

\noindent\textbf{Isomorphism Theorem.} \it Let the boundary-value problem
\eqref{4f1} be parameter-elliptic in the angle $K$. Then there exists a number
$\lambda_{1}\geq1$ that, for each $\lambda\in K$ with $|\lambda|\geq\lambda_{1}$ and
for every $\varphi\in\mathrm{RO}$ subject to \eqref{4f12}, we have the isomorphism
\begin{equation}
(A(\lambda),B(\lambda)):\,H^{\varphi\varrho^{2q}}(\Omega)\,\leftrightarrow\,
\mathcal{H}^{\varphi}(\Omega,\Gamma).
\end{equation}
Moreover, for each fixed $\varphi\in\mathrm{RO}$ satisfying \eqref{4f12}, there
exists a number $c=c(\varphi)\geq1$ such that
\begin{equation}\label{4f15}
\begin{aligned}
c^{-1}\,\|u\|_{\varphi\varrho^{2q},|\lambda|,\Omega}&\leq
\|A(\lambda)u\|_{\varphi,|\lambda|,\Omega}+
\sum_{j=1}^{q}\|B_{j}(\lambda)u\|_{\varphi\varrho^{2q-m_j-1/2},|\lambda|,\Gamma}\\
&\leq c\,\|u\|_{\varphi\varrho^{2q},|\lambda|,\Omega}
\end{aligned}
\end{equation}
for every $\lambda\in K$ with $|\lambda|\geq\lambda_{1}$ and for arbitrary $u\in
H^{\varphi\varrho^{2q}}(\Omega)$. Here, the number $c$ does not depend on $\lambda$
and $u$. \rm

\medskip

In the Sobolev case of $\varphi(t)\equiv t^{s}$, this theorem was proved by
M.~S.~Agranovich and M.~I.~Vishik \cite[\S~4 and 5]{AgranovichVishik64} on the
additional assumption that $s$ is integer. The assumption is not obligatory
\cite[Sec. 3.2]{Agranovich97}. V.~A.~Mikhailets and the second author of the present
paper proved this theorem in the case where the function parameter $\varphi$ varies
regularly in the sense of J.~Karamata at infinity and where each $m_{j}\leq2q-1$
(see \cite[Sec.~7]{07UMJ5} and \cite[Sec.~4.1.4]{MikhailetsMurach10,
MikhailetsMurach14}).

In Section~\ref{sec6}, we will deduce Isomorphism Theorem from the Sobolev case with
the help of the interpolation with a function parameter between Sobolev spaces.

Remark that the right-hand inequality in \eqref{4f15} remains true for every
$\lambda\in\mathbb{C}$ without the assumption that \eqref{4f1} is
parameter--elliptic on $K$ (cf. \cite[Proposition 4.1]{AgranovichVishik64} or
\cite[Sec.~3.2~a]{Agranovich97}).

Isomorphism Theorem implies the following important fact (cf. \cite[Sec.
6.4]{AgranovichVishik64}).

\medskip

\noindent\textbf{Fredholm Property.} \it Suppose that the boundary-value problem
\eqref{4f1} is parameter-elliptic on a certain ray
$K:=\{\lambda\in\mathbb{C}:\arg\lambda =\mathrm{const}\}$. Then the operator
\eqref{4f11} is Fredholm with zero index for each $\lambda\in\mathbb{C}$ and for
every $\varphi\in\mathrm{RO}$ subject to \eqref{4f12}. \rm

\medskip

We will prove these remark and property in Section~\ref{sec6} as well.

Note that the theory of general elliptic boundary-value problems is built for a
narrower class of H\"ormander spaces (called the refined Sobolev scale) by
V.~A.~Mikhailets and A.~A.~Murach in a series of papers, among them we mention the
articles [26--29], survey \cite{12BJMA2}, and monographs \cite{MikhailetsMurach10,
MikhailetsMurach14}. This class consists of the spaces $H^{\varphi}(\Omega)$ such
that $\varphi$ varies regularly at infinity. Specifically, parameter-elliptic
problems is investigated in \cite[Sec.~7]{07UMJ5} under additional assumption that
all $m_{j}\leq2q-1$.

\section{Method of interpolation with a function parameter}\label{sec5}

This method will play a key role in our proof of Isomorphism Theorem. Namely, we
will use the fact that every space $H^{\alpha}(G)$, with $\alpha\in\mathrm{RO}$ and
$G\in\{\mathbb{R}^{n},\Omega,\Gamma\}$, can be obtained by means of the
interpolation (with an appropriate function parameter) between the Sobolev spaces
$H^{(s_0)}(G)$ and $H^{(s_1)}(G)$ with $s_{0}<\sigma_{0}(\alpha)$ and
$s_{1}>\sigma_{1}(\alpha)$.

The method of interpolation with a function parameter between normed spaces was
introduced by C.~Foia\c{s} and J.-L.~Lions in \cite{FoiasLions61} and was then
developed and investigated by some researchers; see monographs
\cite{BrudnyiKrugljak91, Ovchinnikov84} and the literature given therein.

For our purposes it is sufficient to use this interpolation for separable Hilbert
spaces. We recall its definition and some of its properties; see monographs
\cite[Sec 1.1]{MikhailetsMurach10, MikhailetsMurach14} or paper
\cite[Sec.~2]{08MFAT1}, in which this matter is systematically set forth.

Let $X:=[X_{0},X_{1}]$ be a given ordered couple of separable complex Hilbert spaces
$X_{0}$ and $X_{1}$ that satisfy the continuous and dense embedding
$X_{1}\hookrightarrow X_{0}$. We say that this couple is admissible. For $X$ there
exists an isometric isomorphism $J:X_{1}\leftrightarrow X_{0}$ such that $J$ is a
self-adjoint positive operator on $X_{0}$ with the domain $X_{1}$. This operator is
uniquely determined by the couple $X$ and is called a generating operator for~$X$.

By $\mathcal{B}$ we denote the set of all Borel measurable functions
$\psi:(0,\infty)\rightarrow(0,\infty)$ such that $\psi$ is bounded on each compact
interval $[a,b]$, with $0<a<b<\infty$, and that $1/\psi$ is bounded on every set
$[r,\infty)$, with $r>0$.

Let $\psi\in\mathcal{B}$, and consider the operator $\psi(J)$, which is defined (and
positive) in $X_{0}$ as the Borel function $\psi$ of $J$. Denote by
$[X_{0},X_{1}]_{\psi}$ or, simply, by $X_{\psi}$ the domain of the operator
$\psi(J)$ endowed with the inner product
$(u_{1},u_{2})_{X_{\psi}}:=(\psi(J)u_{1},\psi(J)u_{2})_{X_{0}}$ and the
corresponding norm $\|u\|_{X_{\psi}}=\|\psi(J)u\|_{X_{0}}$. The space $X_{\psi}$ is
Hilbert and separable.

A function $\psi\in\mathcal{B}$ is called an interpolation parameter if the
following condition is fulfilled for all admissible couples $X=[X_{0},X_{1}]$ and
$Y=[Y_{0},Y_{1}]$ of Hilbert spaces and for an arbitrary linear mapping $T$ given on
$X_{0}$: if the restriction of $T$ to $X_{j}$ is a bounded operator
$T:X_{j}\rightarrow Y_{j}$ for each $j\in\{0,1\}$, then the restriction of $T$ to
$X_{\psi}$ is also a bounded operator $T:X_{\psi}\rightarrow Y_{\psi}$. In this case
we say that $X_{\psi}$ is obtained by the interpolation with the function parameter
$\psi$ of the couple $X$ (or, otherwise speaking, between the spaces $X_{0}$ and
$X_{1}$).

The function $\psi\in\mathcal{B}$ is an interpolation parameter if and only if
$\psi$ is pseudoconcave on a neighborhood of $+\infty$, i.e.
$\psi(t)\asymp\psi_{1}(t)$ with $t\gg1$ for a certain positive concave
function~$\psi_{1}(t)$. (As usual, the designation $\psi\asymp\psi_{1}$ means that
both the functions $\psi/\psi_{1}$ and $\psi_{1}/\psi$ are bounded on the indicated
set). This fundamental fact follows from J.~Peetre's \cite{Peetre68} description of
all interpolation functions of positive order (see also the monograph
\cite[Sec.~5.4]{BerghLefstrem76}).

Now we can formulate the above-mentioned interpolation property of the extended
Sobolev as follows.

\begin{proposition}\label{prop1}
Let a function $\alpha\in\mathrm{RO}$ and real numbers $s_{0}<\sigma_{0}(\alpha)$
and $s_{1}>\sigma_{1}(\alpha)$ be given. Put
\begin{equation}\label{4f16}
\psi(t):=
\begin{cases}
\;t^{{-s_0}/{(s_1-s_0)}}\,
\alpha\bigl(t^{1/{(s_1-s_0)}}\bigr)&\text{for}\quad t\geq1, \\
\;\alpha(1)&\text{for}\quad0<t<1.
\end{cases}
\end{equation}
Then $\psi\in\mathcal{B}$ is an interpolation parameter, and
\begin{equation}\label{4f17}
\bigl[H^{(s_0)}(G),H^{(s_1)}(G)\bigr]_{\psi}=H^{\alpha}(G)
\end{equation}
with equivalence of norms provided that $G\in\{\mathbb{R}^{n},\Omega,\Gamma\}$.
Moreover, if $G=\mathbb{R}^{n}$, then \eqref{4f17} holds with equality of norms.
\end{proposition}

This proposition is proved in \cite[Theorems 2.18 and 2.22]{MikhailetsMurach10,
MikhailetsMurach14} for $G\in\{\mathbb{R}^{n},\Gamma\}$ and in \cite[Theorem
5.1]{arXiv:1106.2049} for $G=\Omega$.

Besides, the extended Sobolev scale $\{H^{\alpha}(G):\alpha\in\mathrm{RO}\}$
possesses the following important interpolation properties (see \cite[Sec.
2.4.2]{MikhailetsMurach10, MikhailetsMurach14} and \cite[Sec.~2 and 5
]{arXiv:1106.2049}). This scale is closed with respect to the above method of
interpolation with a function parameter and coincides (up to equivalence of norms)
with the class of all Hilbert spaces that are interpolation spaces for the couple of
the Sobolev spaces $H^{(s_0)}(G)$ and $H^{(s_1)}(G)$ with $s_0,s_1\in\mathbb{R}$ and
$s_0<s_1$. The latter property follows from V.~I.~Ovchinnikov's theorem \cite[Sec.
11.4]{Ovchinnikov84}, which gives a description of all Hilbert spaces that are
interpolation spaces for an arbitrarily chosen couple of Hilbert spaces. Recall that
the property of a (Hilbert) space $H$ to be an interpolation space for an admissible
couple $X=[X_0,X_1]$ means the following: a)~the continuous and dense embeddings
$X_1\hookrightarrow H\hookrightarrow X_0$ hold; b)~every linear operator bounded on
each of the spaces $X_0$ and $X_1$ should be bounded on $H$ as well.

Thus, the extended Sobolev scale is the completed extension of the Sobolev scale
$\{H^{(s)}(G):s\in\mathbb{R}\}$ by the interpolation within the category of Hilbert
spaces.

We will also use two properties of interpolation between abstract Hilbert spaces
\cite[Sec 1.1.8 and 1.1.5]{MikhailetsMurach10, MikhailetsMurach14}.

The first of them is an important estimate of the operator norm in interpolation
spaces. Recall that an admissible couple of Hilbert spaces $X=[X_{0},X_{1}]$ is said
to be normal if $\|u\|_{X_{0}}\leq\|u\|_{X_{1}}$ for each $u\in X_{1}$.

\begin{proposition}\label{prop2}
For every interpolation parameter $\psi\in\mathcal{B}$ there exists a number
$\widetilde{c}=\widetilde{c}(\psi)>0$ such that
$$
\|T\|_{X_{\psi}\rightarrow Y_{\psi}}\leq
\widetilde{c}\,\max\,\bigl\{\,\|T\|_{X_{j}\rightarrow Y_{j}}:\,j=0,\,1\bigr\}.
$$
Here, $X=[X_{0},X_{1}]$ and $Y=[Y_{0},Y_{1}]$ are arbitrary normal admissible
couples of Hilbert spaces, and $T$ is an arbitrary linear mapping which is given on
$X_{0}$ and defines the bounded operators $T:X_{j}\rightarrow Y_{j}$, with
$j=0,\,1$. The number $c_{\psi}>0$ does not depend on $X$, $Y$, and $T$.
\end{proposition}

The second property is useful when we interpolate between direct sums of Hilbert
spaces.

\begin{proposition}\label{prop3}
Let $\bigl[X_{0}^{(j)},X_{1}^{(j)}\bigr]$, with $j=1,\ldots,r$, be a finite
collection of admissible couples of Hilbert spaces. Then, for every function
$\psi\in\mathcal{B}$, we have
$$
\biggl[\,\bigoplus_{j=1}^{r}X_{0}^{(j)},\,\bigoplus_{j=1}^{r}X_{1}^{(j)}\biggr]_{\psi}=\,
\bigoplus_{j=1}^{r}\bigl[X_{0}^{(j)},\,X_{1}^{(j)}\bigr]_{\psi}\quad\mbox{with
equality of norms}.
$$
\end{proposition}

\section{Proofs}\label{sec6}

Previously, we will prove an interpolation property of the norm \eqref{4f13} used in
Isomorphism Theorem.

Let $\alpha\in\mathrm{RO}$, with $\sigma_{0}(\alpha)>0$, and let $p\geq1$ and
$G\in\{\mathbb{R}^{n},\Omega,\Gamma\}$. By $H^{\alpha}(G,p)$ we denote the space
$H^{\alpha}(G)$ endowed with the equivalent Hilbert norm \eqref{4f13}. In the
Sobolev case where $\alpha(t)\equiv t^{s}$ for some $s>0$, we also use the notation
$H^{(s)}(G,p)$ and $\|\cdot\|_{(s),p,G}$ for the space $H^{\alpha}(G,p)$ and norm in
it.

\medskip

\noindent\textbf{Interpolation Lemma.} \it Let positive numbers
$s_0<\sigma_{0}(\alpha)$ and $s_1>\sigma_{1}(\alpha)$ be given. Define the
interpolation parameter $\psi$ by formula \eqref{4f16} and denote, for the sake of
brevity, $E(G,p):=\bigl[H^{(s_0)}(G,p),H^{(s_1)}(G,p)\bigr]_{\psi}$. Then there
exists a number $c_{0}\geq1$ such that
\begin{equation}\label{4f18}
c_{0}^{-1}\|w\|_{E(G,p)}\leq\|w\|_{\alpha,p,G}\leq c_{0}\|w\|_{E(G,p)}
\end{equation}
for every $p\geq1$ and all $w\in H^{\alpha}(G)$. The number $c_{0}$ is independent
of $p$ and~$w$. \rm

\medskip

This lemma reinforces Proposition~\ref{prop1} in the sense that the constants in the
estimates of the equivalent norms in $H^{\alpha}(G,p)$ and $E(G,p)$ can be chosen to
be independent of the parameter~$p$.

\begin{proof}[Proof of Interpolation Lemma]
First, we prove this lemma for $G=\mathbb{R}^{n}$. Then we deduce it in the
$G\in\{\Omega,\Gamma\}$ cases, which will be used in our proof of Isomorphism
Theorem.

Suppose that $G=\mathbb{R}^{n}$, with $n\geq1$. Let $p\geq1$ and $w\in
H^{\alpha}(\mathbb{R}^{n})$. Note that
\begin{equation*}
\|w\|_{\alpha,p,\mathbb{R}^{n}}=\biggl(\;\int_{\mathbb{R}^{n}}\,
\bigl(\alpha^{2}(\langle\xi\rangle)+\alpha^{2}(p)\bigr)\,
|\widehat{w}(\xi)|^{2}\,d\xi\,\biggr)^{1/2}.
\end{equation*}
Along with this norm, we consider another Hilbert norm
\begin{equation}\label{4f19}
\|w\|_{\alpha,p,\mathbb{R}^{n}}':=\biggl(\;\int_{\mathbb{R}^{n}}\,
\alpha^{2}(\langle\xi\rangle+p)\,|\widehat{w}(\xi)|^{2}\,d\xi\,\biggr)^{1/2}.
\end{equation}
They are equivalent to each other; moreover,
\begin{equation}\label{4f20}
c_{1}^{-1}\|w\|_{\alpha,p,\mathbb{R}^{n}}'\leq\|w\|_{\alpha,p,\mathbb{R}^{n}}\leq
c_{1}\|w\|_{\alpha,p,\mathbb{R}^{n}}'
\end{equation}
with some number $c_{1}=c_{1}(\alpha)\geq1$ being independent of $p$ and~$w$.

Indeed, since $\sigma_{0}(\alpha)>0$, then there is a number
$c_{2}=c_{2}(\alpha)\geq1$ such that
\begin{equation}\label{4f21}
c_{2}^{-1}\alpha(t_{1}+t_{2})\leq\alpha(t_{1})+\alpha(t_{2})\leq
c_{2}\,\alpha(t_{1}+t_{2})\quad\mbox{for all}\quad t_{1},t_{2}\geq1.
\end{equation}
This immediately implies \eqref{4f20} with $c_{1}:=\sqrt{2}\,c_{2}$. The
demonstration of property \eqref{4f21} is simple and will be given at once after the
proof of this lemma.

Let $H^{\alpha}(\mathbb{R}^{n},p,1)$ denote the space $H^{\alpha}(\mathbb{R}^{n})$
endowed with the equivalent Hilbert norm \eqref{4f19}. In the Sobolev case where
$\alpha(t)\equiv t^{s}$ for some $s>0$, we also use the notation
$H^{(s)}(\mathbb{R}^{n},p,1)$ and $\|\cdot\|_{(s),p,\mathbb{R}^{n}}'$ for the space
$H^{\alpha}(\mathbb{R}^{n},p,1)$ and norm in it. Put, for the sake of brevity,
$$
E(\mathbb{R}^{n},p,1):=
\bigl[H^{(s_0)}(\mathbb{R}^{n},p,1),H^{(s_1)}(\mathbb{R}^{n},p,1)\bigr]_{\psi}.
$$

Let us prove the equality
\begin{equation}\label{4f22}
\|w\|_{E(\mathbb{R}^{n},p,1)}=\|w\|_{\alpha,p,\mathbb{R}^{n}}'.
\end{equation}

To this end, calculate the norm in the interpolation space $E(\mathbb{R}^{n},p,1)$.
It is directly verified that the pseudodifferential operator
$$
J:\,u\mapsto\mathcal{F}^{-1}
\bigl[(\langle\xi\rangle+p)^{s_1-s_0}\,\widehat{u}(\xi)\bigr],\quad\mbox{with}\quad
u\in H^{(s_1)}(\mathbb{R}^{n}),
$$
is the generating operator for the admissible couple of the Hilbert spaces
$H^{(s_0)}(\mathbb{R}^{n},p,1)$ and $H^{(s_1)}(\mathbb{R}^{n},p,1)$. Here,
$\mathcal{F}^{-1}$ is the inverse Fourier transform in $\xi\in\mathbb{R}^{n}$.
Applying the isometric isomorphism
$$
\mathcal{F}:H^{(s_{0})}(\mathbb{R}^{n},p,1)\leftrightarrow
L_{2}\bigl(\mathbb{R}^{n},\,(\langle\xi\rangle+p)^{2s_{0}}\,d\xi\bigr),
$$
we reduce the operator $\psi(J)$ to the form of multiplication by the function
$$\psi\bigl((\langle\xi\rangle+p)^{s_1-s_0}\bigr)=
(\langle\xi\rangle+p)^{-s_0}\,\alpha(\langle\xi\rangle+p)\quad\mbox{of}\quad\xi,
$$
in view of \eqref{4f16}. Therefore,
\begin{gather*}
\|w\|_{E(\mathbb{R}^{n},p,1)}=
\|\psi(J)w\|_{(s_{0}),p,\mathbb{R}^{n}}'= \\
\biggl(\,\int_{\mathbb{R}^{n}}\psi^{2}\bigl((\langle\xi\rangle+p)^{s_1-s_0}\bigr)\,
|\widehat{w}(\xi)|^{2}\,(\langle\xi\rangle+p)^{2s_{0}}\,d\xi\biggr)^{1/2}= \\
\biggl(\,\int_{\mathbb{R}^{n}}
\alpha^{2}(\langle\xi\rangle+p)\,|\widehat{w}(\xi)|^{2}\,d\xi\biggr)^{1/2}=
\|w\|_{\alpha,p,\mathbb{R}^{n}}'.
\end{gather*}
Equality \eqref{4f22} is proved.

We note that
\begin{equation}\label{4f23}
c_{3}^{-1}\,\|w\| _{E(\mathbb{R}^{n},p)}\leq \|w\|_{E(\mathbb{R}^{n},p,1)}\leq
c_{3}\,\|w\|_{E(\mathbb{R}^{n},p)}
\end{equation}
with some number $c_{3}\geq1$, this number being independent of $p$ and $w$. Indeed,
the identity operator $I$ sets the isomorphism
\begin{equation*}
I:\,H^{(s)}(\mathbb{R}^{n},p)\leftrightarrow
H^{(s)}(\mathbb{R}^{n},p,1)\quad\mbox{for each}\quad s>0.
\end{equation*}
The norms of this operator and its inverse do not exceed $2^{s}$. Consider this
isomorphism for each $s\in\{s_{0},s_{1}\}$ and apply the interpolation with the
function parameter $\psi$. We obtain the isomorphism
\begin{gather*}
I:\,E(\mathbb{R}^{n},p)\leftrightarrow E(\mathbb{R}^{n},p,1)
\end{gather*}
and conclude by Proposition~\ref{prop2} that the norms of this operator and its
inverse do not exceed $2^{s_{1}}\widetilde{c}$, where the number
$\widetilde{c}=\widetilde{c}(\psi)$ depends only on $\psi$. This means \eqref{4f23}
with $c_{3}:=2^{s_{1}}\widetilde{c}$.

Formulas \eqref{4f20}, \eqref{4f22}, and \eqref{4f23} directly imply the required
estimate \eqref{4f18} for $G=\mathbb{R}^{n}$ and $c_{0}:=c_{1}c_{3}$, with $c_{0}$
being independent of $p$ and $w$. Thus,
\begin{equation}\label{4f24}
(c_{1}c_{3})^{-1}\|w\|_{E(\mathbb{R}^{n},p)}\leq\|w\|_{\alpha,p,\mathbb{R}^{n}}\leq
c_{1}c_{3}\|w\|_{E(\mathbb{R}^{n},p)}
\end{equation}
for all $w\in H^{\alpha}(\mathbb{R}^{n})$ and $p\geq1$.

Now, let us deduce Interpolation Lemma in the case of $G=\Omega$. As above,
$p\geq1$. Let $R$ denote the operator that restricts distributions
$w\in\mathcal{S}'(\mathbb{R}^{n})$ to $\Omega$. We have the bounded linear operators
\begin{gather} \label{4f25}
R:\,H^{\alpha}(\mathbb{R}^{n},p)\rightarrow
H^{\alpha}(\Omega,p), \\
R:\,H^{(s)}(\mathbb{R}^{n},p)\rightarrow H^{(s)}(\Omega,p),\quad\mbox{with}\quad
s>0. \label{4f26}
\end{gather}
Certainly, their norms do not exceed $1$.

Consider the operator \eqref{4f26} for each $s\in\{s_{0},s_{1}\}$ and apply the
interpolation with the parameter $\psi$. We obtain the bounded operator
\begin{equation*}
R:\,E(\mathbb{R}^{n},p)\rightarrow E(\Omega,p), \quad\mbox{whose
norm}\;\leq\widetilde{c}
\end{equation*}
by Proposition~\ref{prop2}. Hence, in view of \eqref{4f24}, we have the bounded
operator
\begin{equation}\label{4f27}
R:\,H^{\alpha}(\mathbb{R}^{n},p)\rightarrow E(\Omega,p), \quad\mbox{whose
norm}\;\leq c_{4}:=c_{1}c_{3}\widetilde{c}.
\end{equation}

Further, we need to use a bounded linear operator which is a right-inverse of both
the operators \eqref{4f25} and \eqref{4f27}. As is known, for each integer $k\geq1$
there exists a linear mapping $T_{k}$ that extends every function $u\in
H^{(0)}(\Omega)$ onto $\mathbb{R}^{n}$ and that sets the bounded operator
\begin{equation}\label{4f28}
T_{k}:\,H^{(s)}(\Omega)\rightarrow H^{(s)}(\mathbb{R}^{n})\quad\mbox{for every}\quad
s\in[0,k].
\end{equation}
Such a mapping is constructed in \cite[Sec~4.2.2]{Triebel95}, for example.

Let $k:=[s_{1}]+1$, and consider the operator \eqref{4f28} for each
$s\in\{s_{0},s_{1}\}$. Applying the interpolation with the parameter $\psi$ and
Proposition~\ref{prop1}, we obtain the bounded operator
\begin{equation}\label{4f29}
T_{k}:\,H^{\alpha}(\Omega)\rightarrow H^{\alpha}(\mathbb{R}^{n}).
\end{equation}
Let $c_{5}$ stand for the maximum of the norms of operators \eqref{4f29} and
\eqref{4f28} with $s\in\{0,s_{0},s_{1}\}$. We have the bounded operators
\begin{equation}\label{4f30}
T_{k}:\,H^{\alpha}(\Omega,p)\rightarrow H^{\alpha}(\mathbb{R}^{n},p),
\quad\mbox{whose norm}\;\leq c_{5},
\end{equation}
and
\begin{equation}\label{4f31}
T_{k}:\,H^{(s)}(\Omega,p)\rightarrow H^{(s)}(\mathbb{R}^{n},p),\quad
s\in\{s_{0},s_{1}\},\quad\mbox{whose norms}\;\leq c_{5}.
\end{equation}
Applying the interpolation with the parameter $\psi$ and Proposition~\ref{prop2} to
\eqref{4f31}, we obtain the bounded operator
$$
T_{k}:\,E(\Omega,p)\rightarrow E(\mathbb{R}^{n},p),\quad\mbox{whose norm}\;\leq
c_{5}\widetilde{c}.
$$
Hence, in view of \eqref{4f24}, we have the bounded operator
\begin{equation}\label{4f32}
T_{k}:\,E(\Omega,p)\rightarrow H^{\alpha}(\mathbb{R}^{n},p),\quad\mbox{whose
norm}\;\leq c_{6}:=c_{1}c_{3}c_{5}\widetilde{c}.
\end{equation}

Now, considering the product of the operators \eqref{4f32} and \eqref{4f25}, and
then, the product of the operators \eqref{4f30} and \eqref{4f27}, we arrive at the
identity operators
$$
I=RT_{k}:\,E(\Omega,p)\rightarrow H^{\alpha}(\Omega,p),\quad\mbox{whose norm}\;\leq
c_{6},
$$
and
$$
I=RT_{k}:\,H^{\alpha}(\Omega,p)\rightarrow E(\Omega,p),\quad\mbox{whose norm}\;\leq
c_{4}c_{5}.
$$
This gives the required estimate \eqref{4f18} for $G=\Omega$ and
$c_{0}:=\max\{c_{4}c_{5},c_{6}\}$, with $c_{0}$ being independent of $p$ and~$w$.

Finally, let us deduce Interpolation Lemma in the case where $G=\Gamma$. We fix a
finite collection of local charts $\{\alpha_{j}\}$ on $\Gamma$ and corresponding
partition of unity $\{\chi_{j}\}$ used in the definition of the extended Sobolev
scale on~$\Gamma$; here, $j$ runs over the values $1,\ldots,\varkappa$.

Consider the linear mapping
$$
T:\,h\mapsto\bigl((\chi_{1}h)\circ\alpha_{1},\ldots,
(\chi_{\varkappa}h)\circ\alpha_{\varkappa}\bigr),
$$
with $h\in\mathcal{D}'(\Gamma)$. We can directly verify that this mapping defines
the isometric operators
\begin{gather}\label{4f33}
T:\,H^{\alpha}(\Gamma,p)\rightarrow
\bigl(H^{\alpha}(\mathbb{R}^{n-1},p)\bigr)^{\varkappa}, \\
T:\,H^{(s)}(\Gamma,p)\rightarrow
\bigl(H^{(s)}(\mathbb{R}^{n-1},p)\bigr)^{\varkappa}, \quad\mbox{with}\quad s>0.
\label{4f34}
\end{gather}
Consider the operator \eqref{4f34} for each $s\in\{s_{0},s_{1}\}$ and apply the
interpolation with the parameter $\psi$. According to Propositions \ref{prop2} and
\ref{prop3}, we obtain the bounded operator
\begin{equation*}
T:\,E(\Gamma,p)\rightarrow
\bigl(E(\mathbb{R}^{n-1},p)\bigr)^{\varkappa},\quad\mbox{whose norm}\;\leq
\widetilde{c}.
\end{equation*}
Hence, in view of \eqref{4f24}, we have the bounded operator
\begin{equation}\label{4f35}
T:\,E(\Gamma,p)\rightarrow \bigl(H^{\alpha}(\mathbb{R}^{n-1},p)\bigr)^{\varkappa},
\quad\mbox{whose norm}\;\leq c_{4}=c_{1}c_{3}\widetilde{c}.
\end{equation}

Along with $T$, consider the linear mapping
$$
Q:\,(w_{1},\ldots,w_{\varkappa})\mapsto
\sum_{j=1}^{\varkappa}\Theta_{j}\bigl((\eta_{j}w_{j})\circ\alpha_{j}^{-1}\bigr),
$$
with $w_{1},\ldots,w_{\varkappa}\in\mathcal{S}'(\mathbb{R}^{n})$. Here, the function
$\eta_{j}\in C^{\infty}(\mathbb{R}^{n})$ is equal to $1$ on the set
$\alpha_{j}^{-1}(\mathrm{supp}\,\chi_{j})$ and is compactly supported, whereas
$\Theta_{j}$ denotes the operator of extension by zero from $\Gamma_{j}$ onto
$\Gamma$. The mapping $Q$ is a left-inverse of $T$. Indeed,
$$
QTh=\sum_{j=1}^{\varkappa}\,\Theta_{j}\bigl((\eta_{j}\,
((\chi_{j}h)\circ\alpha_{j}))\circ\alpha_{j}^{-1}\bigr)= \sum_{j=1}^{\varkappa}\,
\Theta_{j}\bigl((\chi_{j}h)\circ\alpha_{j}\circ\alpha_{j}^{-1}\bigr)=
\sum_{j=1}^{\varkappa}\,\chi_{j}h=h
$$
for every $h\in\mathcal{D}'(\Gamma)$.

We have the bounded operators
\begin{gather}\label{4f36}
Q:\,\bigl(H^{\alpha}(\mathbb{R}^{n-1})\bigr)^{\varkappa}\rightarrow
H^{\alpha}(\Gamma),\\
Q:\,\bigl(H^{(s)}(\mathbb{R}^{n-1})\bigr)^{\varkappa}\rightarrow
H^{(s)}(\Gamma),\quad\mbox{with}\quad s\in\mathbb{R}. \label{4f37}
\end{gather}
The boundedness of \eqref{4f37} is a known property of Sobolev spaces (see, e.g.,
\cite[Sec. 2.6]{Hermander63}). The boundedness of \eqref{4f36} follows from this
property with the help of interpolation by virtue of Propositions \ref{prop1}
and~\ref{prop3}.

Let $c_{7}$ denote the maximum of the norms of operators \eqref{4f36} and
\eqref{4f37} with $s\in\{0,s_{0},s_{1}\}$. We have the bounded operators
\begin{equation}\label{4f38}
Q:\,\bigl(H^{\alpha}(\mathbb{R}^{n-1},p)\bigr)^{\varkappa}\rightarrow
H^{\alpha}(\Gamma,p), \quad\mbox{whose norm}\;\leq c_{7},
\end{equation}
and
\begin{equation}\label{4f39}
Q:\,\bigl(H^{(s)}(\mathbb{R}^{n-1},p)\bigr)^{\varkappa}\rightarrow
H^{(s)}(\Gamma,p),\quad s\in\{s_{0},s_{1}\},\quad\mbox{whose norms}\;\leq c_{7}.
\end{equation}
Applying the interpolation with the parameter $\psi$ to \eqref{4f39}, we obtain the
bounded operator
\begin{equation*}
Q:\,\bigr(E(\mathbb{R}^{n-1},p)\bigr)^{\varkappa}\rightarrow
E(\Gamma,p),\quad\mbox{whose norm}\;\leq c_{7}\widetilde{c}
\end{equation*}
by virtue of Propositions \ref{prop2} and~\ref{prop3}. Hence, in view of
\eqref{4f24}, we have the bounded operator
\begin{equation}\label{4f40}
Q:\,\bigl(H^{\alpha}(\mathbb{R}^{n-1},p)\bigr)^{\varkappa}\rightarrow E(\Gamma,p),
\quad\mbox{whose norm}\;\leq c_{8}:=c_{1}c_{3}c_{7}\widetilde{c}.
\end{equation}

Thus, multiplying \eqref{4f40} by the isometric operator \eqref{4f33}, we obtain the
bounded identity operator
$$
I=QT:\,H^{\alpha}(\Gamma,p)\rightarrow E(\Gamma,p),\quad\mbox{whose norm}\;\leq
c_{8}.
$$
Besides, taking the product of the operators \eqref{4f38} and \eqref{4f35}, we get
the bounded identity operator
$$
I=QT:\,E(\Gamma,p)\rightarrow H^{\alpha}(\Gamma,p),\quad\mbox{whose norm}\;\leq
c_{4}c_{7}.
$$
This yields the required estimate \eqref{4f18} for $G=\Gamma$ and
$c_{0}:=\max\{c_{8},c_{4}c_{7}\}$, with $c_{0}$ being independent of $p$ and~$w$.
\end{proof}

As we have promised, let us prove that every function $\alpha\in\mathrm{RO}$ with
$\sigma_{0}(\alpha)>0$ satisfies property \eqref{4f21}.

Since $\sigma_{0}(\alpha)>0$, then, by \eqref{4f6}, there exists a number $c'>0$
such that
\begin{equation}\label{4f41}
c'\leq\frac{\alpha(\lambda t)}{\alpha(t)}\quad\mbox{for all}\quad
t\geq1\quad\mbox{and}\quad\lambda\geq1.
\end{equation}
Hence, for arbitrary $t_{1},t_{2}\geq1$, we can write
$c'\alpha(t_{j})\leq\alpha(t_{1}+t_{2})$ with $j\in\{1,\,2\}$, which gives the
right-hand inequality in \eqref{4f21}.

Besides, applying properties \eqref{4f41} and \eqref{4f5} of $\alpha\in\mathrm{RO}$,
we deduce the left-hand inequality in \eqref{4f21}, namely,
$$
c'\alpha(t_{1}+t_{2})\leq\alpha(2t)\leq c''\alpha(t)\leq
c''(\alpha(t_{1})+\alpha(t_{2}))\quad\mbox{for all}\quad t_{1},t_{2}\geq1.
$$
Here, $t:=\max\{t_{1},t_{2}\}$, and $c''$ is a certain positive number, which does
not depend on $t_{1}$ and $t_{2}.$ Thus, property \eqref{4f21} is proved.

Now, applying Interpolation Lemma, we will give

\begin{proof}[Proof of Isomorphism Theorem]
In the Sobolev case of $\varphi(t)\equiv t^{s}$, this theorem was proved by
M.~S.~Agranovich and M.~I.~Vishik (see \cite[\S~4 and 5]{AgranovichVishik64} and
\cite[Sec. 3.2]{Agranovich97}). Using spaces $H^{(s)}(G,|\lambda|)$, we can
reformulate their result in the following way. Let the boundary-value problem
\eqref{4f1} be parameter-elliptic in the angle $K$. Then there exists a number
$\lambda_{1}\geq1$ that the isomorphism
\begin{align}\label{4f42}
(A(\lambda),B(\lambda))&:\,H^{(s+2q)}(\Omega,|\lambda|)\leftrightarrow
H^{(s)}(\Omega,|\lambda|)\oplus
\bigoplus_{j=1}^{q}H^{(s+2q-m_j-1/2)}(\Gamma,|\lambda|)\\
&=:\mathcal{H}^{(s)}(\Omega,\Gamma,|\lambda|)\notag
\end{align}
holds for every $s>l$ and each $\lambda\in K$ with $|\lambda|\geq\lambda_{1}$.
Moreover, the norms of the operator \eqref{4f42} and its inverse are uniformly
bounded with respect to $\lambda$.

Deduce Isomorphism Theorem for the extended Sobolev scale. Let
$\varphi\in\mathrm{RO}$ meet condition \eqref{4f12}, and let $\lambda\in K$ satisfy
$|\lambda|\geq\lambda_{1}$. We choose $l_{0},l_{1}\in\mathbb{R}$ so that
$l<l_{0}<\sigma_{0}(\varphi)$ and $l_{1}>\sigma_{1}(\varphi)$. Define the
interpolation parameter $\psi\in\mathcal{B}$ by formula \eqref{4f16}, in which
$\alpha:=\varphi$, $s_{0}:=l_{0}$, and $s_{1}:=l_{1}$.

Consider the isomorphism \eqref{4f42} for each $s\in\{l_{0},l_{1}\}$ and apply the
interpolation with the function parameter $\psi$. According to Propositions
\ref{prop2} and \ref{prop3}, we obtain the isomorphism
\begin{align}\notag
(A(\lambda),B(\lambda))&:\, \bigl[H^{(l_{0}+2q)}(\Omega,|\lambda|),
H^{(l_{1}+2q)}(\Omega,|\lambda|)\bigr]_{\psi}\leftrightarrow
\bigl[H^{(l_{0})}(\Omega,|\lambda|),H^{(l_{1})}(\Omega,|\lambda|)\bigr]_{\psi}\\
&\oplus\bigoplus_{j=1}^{q}\bigl[H^{(l_{0}+2q-m_{j}-1/2)}(\Gamma,|\lambda|),
H^{(l_{1}+2q-m_{j}-1/2)}(\Gamma,|\lambda|)\bigr]_{\psi}.\label{4f43}
\end{align}
The norms of the operator \eqref{4f43} and its inverse are uniformly bounded with
respect to~$\lambda$.

Hence, we deduce by Interpolation Lemma that the isomorphism \eqref{4f43} acts
between the following spaces:
\begin{align}\label{4f44}
(A(\lambda),B(\lambda))&:\,H^{\varphi\varrho^{2q}}(\Omega,|\lambda|)\leftrightarrow
H^{\varphi}(\Omega,|\lambda|)\oplus
\bigoplus_{j=1}^{q}H^{\varphi\varrho^{2q-m_j-1/2}}(\Gamma,|\lambda|)\\
&:=\mathcal{H}^{\varphi}(\Omega,\Gamma,|\lambda|).\notag
\end{align}
Moreover, the norms of the operator \eqref{4f44} and its inverse are uniformly
bounded with respect to~$\lambda$. Note that we have also applied this lemma in the
case where
$$
\alpha=\varphi\varrho^{2q},\quad s_{1}=l_{1}+2q,\quad\mbox{and}\quad s_{0}=l_{0}+2q,
$$
and in the case where
$$
\alpha=\varphi\varrho^{2q-m_{j}-1/2},\quad
s_{0}=l_{0}+2q-m_{j}-1/2,\quad\mbox{and}\quad s_{1}=l_{1}+2q-m_{j}-1/2.
$$
In these cases the function $\psi$ still satisfies \eqref{4f16}.

Owing to the properties of the operator \eqref{4f44}, we can state that Isomorphism
Theorem is proved.
\end{proof}

Analysing this proof, we note the following. According to
\cite[Sec.~3.2~a]{Agranovich97}, the mapping $u\mapsto(A(\lambda)u,B(\lambda)u)$,
with $u\in C^\infty(\overline{\Omega})$, extends uniquely (by continuity) to the
bounded linear operator
\begin{equation}\label{4f45}
(A(\lambda),B(\lambda)):\,H^{(s+2q)}(\Omega,|\lambda|)\rightarrow
\mathcal{H}^{(s)}(\Omega,\Gamma,|\lambda|)
\end{equation}
for every $s>l$ and each $\lambda\in\mathbb{C}$. Moreover, the norm of this operator
is uniformly bounded with respect to~$\lambda$. Here, we do not suppose that the
boundary-value problem \eqref{4f1} is parameter-elliptic in $K$. Applying the
interpolation and reasoning as in the latter proof, we obtain the bounded operator
\begin{equation}\label{4f46}
(A(\lambda),B(\lambda)):\,H^{\varphi\varrho^{2q}}(\Omega,|\lambda|)\rightarrow
\mathcal{H}^{\varphi}(\Omega,\Gamma,|\lambda|)
\end{equation}
for each $\lambda\in\mathbb{C}$ and every $\varphi\in\mathrm{RO}$ with
$\sigma_{0}(\varphi)>l$. Its norm is uniformly bounded with respect to~$\lambda$.
This means that the right-hand inequality in \eqref{4f15} remains true for indicated
$\lambda$ and $\varphi$ without the assumption about parameter--ellipticity.

It remains to give

\begin{proof}[Proof of Fredholm Property]
Let the boundary-value problem \eqref{4f1} be parameter-elliptic on a certain closed
ray $K:=\{\lambda\in\mathbb{C}:\arg\lambda=\mathrm{const}\}$. According to
\cite[Sec. 6.4]{AgranovichVishik64}, the bounded operator \eqref{4f45} is Fredholm
with zero index for every $s>l$ and each $\lambda\in\mathbb{C}$. Moreover, the
kernel of this operator lies in $C^{\infty}(\overline{\Omega})$ and, hence, does not
depend on~$s$. Both the Fredholm property and index are preserved under
interpolation of normed spaces and bounded operators provided that these operators
have the same index and the same kernel (see, e.g.,
\cite[Sec.~1.1.7]{MikhailetsMurach10, MikhailetsMurach14}). Therefore, using the
interpolation with a function parameter and applying Proposition~\ref{prop1}, we
conclude that the bounded operator \eqref{4f46} is Fredholm with zero index for each
$\lambda\in\mathbb{C}$ and every $\varphi\in\mathrm{RO}$ with
$\sigma_{0}(\varphi)>l$.
\end{proof}

\end{document}